\documentclass[11pt]{article}
\usepackage{graphicx}
\usepackage{amsmath,amssymb,amsthm,amsfonts}
\usepackage{amssymb}
\newtheorem{thm}{Theorem}[section]

\newtheorem{lem}[thm]{Lemma}

\theoremstyle{definition}

\theoremstyle{remark}

\numberwithin{equation}{section}

\DeclareMathSymbol{\C}{\mathalpha}{AMSb}{"43}

\newcommand{\eps}{\varepsilon}

\newcommand{\lam}{\lambda}
\newcommand{\alp}{\alpha}

\newcommand{\R}{{\mathbb{R}}}

\def\endproof{\hfill$\blacksquare$\vspace{6pt}}

\newcommand{\bsub}{\begin{subequations}}
\newcommand{\esub}{\end{subequations}$\!$}
\begin{document}

\title{Estimates for the quenching time of a parabolic
equation modeling electrostatic MEMS}\author{ Nassif
Ghoussoub\thanks{Partially supported by the Natural Science
and Engineering Research Council of Canada.} \\ Department of Mathematics, University of British Columbia,\\
Vancouver, B.C. Canada V6T 1Z2\\ [3mm] Yujin
Guo  \\
School of Mathematics, University of Minnesota,\\
 Minneapolis, MN 55455 USA
}
\date{}
\smallbreak \maketitle

\begin{abstract}
\noindent The singular parabolic problem  $u_t=\Delta u
-\frac{\lambda f(x)}{(1+u)^2}$ on a bounded domain $\Omega$ of
$\R^N$ with Dirichlet boundary conditions, models the dynamic
deflection of an elastic  membrane in a simple electrostatic
Micro-Electromechanical System (MEMS) device. In this paper,  we
analyze and estimate the quenching time of the elastic membrane in
terms of the applied voltage ---represented here by $\lam$. As a
byproduct,  we prove               that for sufficiently large
$\lam$, finite-time quenching  must occur near the maximum point of
the varying dielectric permittivity profile $f(x)$.
\end{abstract}

\vskip 0.2truein


\section{Introduction}

\begin{figure}[htbp]
\begin{center}
\includegraphics*[width = 12cm, height = 5cm, clip]{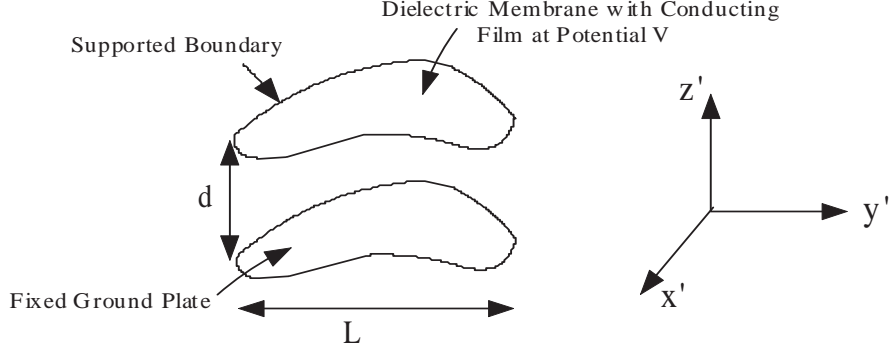}
\caption{{\em The simple electrostatic MEMS device.}}
\label{fig:fig1}
\end{center}
\end{figure}

Micro-Electromechanical Systems (MEMS) are often used to combine
electronics with micro-size mechanical devices in the design of
various types of microscopic machinery.
An overview of the physical phenomena of the mathematical models
associated with the rapidly developing field of MEMS technology is
given in \cite{PB}. The key component of many modern MEMS is the
simple idealized electrostatic device shown in
Figure~\ref{fig:fig1}. The upper part of this device consists of a
thin and deformable elastic membrane that is held fixed along its
boundary and which lies above a rigid grounded plate. This elastic
membrane is modeled as a dielectric with a small but finite
thickness. The upper surface of the membrane is coated with a
negligibly thin metallic conducting film. When a voltage $V$ is
applied to the conducting film, the thin dielectric membrane
deflects towards the bottom plate, and when $V$ is increased beyond
a certain critical value $V^*$ --known as pull-in voltage--  the
steady-state of the elastic membrane is lost, and proceeds to
quenching, $i.e.$ snap through, at a finite time creating the
so-called pull-in instability.

   A mathematical model of the physical phenomena, leading to a partial
differential equation for the dimensionless dynamic deflection of
the membrane, was derived and analyzed in  \cite{GPW}. In the
damping-dominated limit, and using a narrow-gap asymptotic analysis,
the dimensionless dynamic deflection $u=u(x,t)$ of the membrane on a
bounded domain $\Omega$ in $\mathbb{R}^2$, is found to satisfy the
following parabolic problem   $$  \arraycolsep=1.5pt
\arraycolsep=1.5pt \begin {array}{lll} u_t-\Delta u& =\displaystyle
\frac{\lambda f(x)}{(1-u)^2} \quad& {\rm for}
\  x \in \Omega \,,\\[3mm]
\hfill u (x,t) &=0  \quad  & {\rm for}
\ x \in \partial\Omega \,,\\
\hfill u(x,0) &=0  \quad &{\rm for} \  x \in \Omega.
\end{array}
 \eqno{(P)_\lam} $$   
The initial condition in $(P)_\lam$
assumes that the membrane is initially undeflected and the voltage
is suddenly applied to the upper surface of the membrane at time
$t=0$. The parameter $\lambda >0$ in $(P)_\lam$ characterizes the
relative strength of the electrostatic and mechanical forces in the
system, and is given in terms of the applied voltage $V$ by $\lambda
=\frac{\eps_0 V^2 L^2}{2T_e d^3}$, where $d$ is the undeflected gap
size, $L$ is the length scale of the membrane, $T_e$ is the tension
of the membrane, and $\eps_0$ is the permittivity of free space in
the gap between the membrane and the bottom plate. We shall  use
from now on the parameter $\lam $ and $\lam ^*$ to represent the
applied voltage $V$ and pull-in voltage $V^*$, respectively.
Referred to as the {\em permittivity profile}, $f(x)$ in $(P)_\lam$ is
defined by the ratio $   f(x) = \frac{\eps_0}{\eps_{2}(x)}$, where
$\eps_{2}(x)$ is the dielectric permittivity of the thin membrane.

Consider first the steady-state solutions
 of $(P)_\lam$
$$
\arraycolsep=1.5pt\begin{array}{lll}\arraycolsep=1.5pt
   \hfill   -\Delta w &=& \displaystyle\frac{\lambda
f(x)}{(1-w)^2} \,   \  \ x \in \Omega, \\
 \hfill w (x)  &=& 0  \quad\quad\quad\quad x \in \partial\Omega
\end{array}\eqno{(S)_{\lam}}$$
with $0<w<1$ on $\Omega \subset\mathbb{R}^N$, and $f(x)$ was assumed to satisfy
\begin{equation}\label{IV:permittivity}
\begin{array}{lll}\arraycolsep=1.5pt
\hbox{$f\in C^\alpha(\bar \Omega)$ for some $\alpha \in (0,1]$,
$0\leq f\leq 1$  and }\\\hbox{$f>0$ on a subset of $\Omega$ with
positive measure.}
\end{array}
\end{equation}
One can then easily show (e.g., Theorem 1.1 in \cite{GG1}) that there exists a finite pull-in voltage $\lambda^*:=\lam
^*(\Omega, f)>0$ such that:
\begin{trivlist}
\item $\bullet$\,  If $0\leq \lam <\lam ^*$, there exists at least one solution for
$(S)_{ \lam}$.
\item $\bullet$\, If  $\lam >\lam ^*$, there is no solution for $(S)_{\lam}$.
 \end{trivlist}


 \noindent Upper and lower  bounds on the pull-in voltage $\lam ^*$ were also
given in Theorem 1.1 of \cite{GG1}. Fine properties of steady states
--such as regularity, stability, uniqueness, multiplicity, energy
estimates and comparison results-- were shown in \cite{EGG} and
\cite{GG1} to depend on the dimension of the ambient space and on
the permittivity profile.

For the dynamic problem $(P)_\lam$, we first define the following notion.

\noindent {{\bf Definition 1.1.}\   \em (1)\, A solution $u(x,t)$ of
$(P)_\lam$ is said to be quenching  at a --possibly infinite-- time
$T=T(\lam , f, \Omega  )$,  if the maximal value of $u$ reaches $1$
at  time $T$.

(2)\,  A point $x_0\in \bar \Omega$ is said to be a {\it quenching
point} for a solution $u(x,t)$ of $(P)_\lam$, if for some $T\in (0,
+\infty]$, we have $\lim\limits_{t_n\to T} u(x_0,t_n)=1$. }

\noindent In \cite{GG2} we dealt with issues of global convergence
as well as quenching in finite or infinite time  of the solutions of
$(P)_\lam$. One of the main results was the following relationship
between the voltage $\lam$ and the nature of the dynamic solution
$u$ of $(P)_\lam$.

 \noindent {\bf Theorem A (Theorem 1.1 in \cite{GG2}).}
{\em Assuming that  $f$ satisfies $(1.1)$ on a bounded domain
$\Omega$, then the
followings hold:
  \begin{enumerate}
  \item  If $\lam \le\lam ^*$, then there exists a unique
solution $u(x,t)$ for $(P)_\lam$ which globally converges pointwise as
$t\to +\infty $ to its unique minimal steady-state.
\item  If $\lam >\lam ^*$  and $\inf_\Omega f>0$,  then the  unique solution  $u(x,t)$ of $(P)_\lam$ must be quenching at a finite time.
   \end{enumerate}}
 \noindent A refined description of finite-time quenching behavior for $u$ was given in \cite{GG3}, where some quenching estimates, quenching rates, as well as some
information on the properties of quenching set --such as
compactness, location and shape, were obtained.

The first purpose of this paper is to prove --in Theorem 2.1-- that
quenching in  finite-time occurs as soon as $\lam >\lam ^*$, which
means that Theorem A. 2. above holds without the restriction
$\inf_\Omega f>0$. On the other hand, we continue our search for
optimal estimates on quenching times at voltages $\lam >\lam ^*$,
since the latter translate into useful information on the operation
speed of MEMS devices. Indeed, we established in Theorem 1.3 of
\cite{GG2}, that if $\inf_{x\in \Omega} f(x)>0$, then  the following
upper estimate for the quenching time holds  for any
$\lambda>\lambda^*$:
\begin{equation}\label{up.estimate.0}
T_\lambda(\Omega, f)\le 
\frac{8(\lam +\lam ^*)^2}{3\inf_{x\in \Omega} f(x)(\lam -\lam ^*)^2(\lam
+3\lam ^*)} \Big[1+\Big(\frac{\lam +3\lam ^*}{2\lam
+2\lam^*}\Big)^{1/2}\Big].
\end{equation}
In this paper, we shall improve this estimate --at least in dimensions less than 8-- by proving that
$$ T_\lambda(\Omega, f)\sim C\big(\lam -\lam
^*\big)^{-\frac{1}{2}}\quad \quad \mbox{as} \quad \lam \searrow\lam
^*,$$ while
$$
T\sim \frac{1}{3\lam \sup_{x\in \bar \Omega }f(x)} \quad \mbox{as}
\quad \lam \nearrow\infty .
$$

To be more precise, we first recall the following notions and results from \cite{GG1}.
For any solution $w$ of $(S)_{\lam}$, we consider the linearized operator at $w$ defined by
$L_{w, \lambda} =-\Delta -\frac{2\lambda f(x)}{(1-w)^3},$ and its
corresponding eigenvalues $\{\mu_{k, \lambda}(w); k=1, 2,...\}$. Say that a solution
$w_{\lam}$ of $(S)_{ \lam}$ is minimal, if
 $w_{\lam}(x)\le w(x)$ in $\Omega $ whenever $w$ is any solution of
 $(S)_{\lam}$. We recall the following


\noindent {\bf Theorem B (Theorem 1.2 in \cite{GG1}).} {\em Assume
$f$ satisfies $(\ref{IV:permittivity})$ on a bounded domain $\Omega
\subset \mathbb{R}^N$.
Then,
\begin{trivlist}
\item 1.\, For any $0\leq \lam <\lam ^*$, there exists a unique minimal solution $w_{\lam}$
of $(S)_{ \lam}$ such that $\mu_{1,\lambda}(w_\lambda)>0$. Moreover
for each $x\in \Omega$, the function $\lambda \to w_\lambda (x)$ is
strictly increasing and differentiable on $(0, \lambda^*)$.
\item 2.\, If $1\leq N\leq 7$,  then
$w^*=\displaystyle\lim_{\lambda \uparrow
\lambda^*}w_\lambda$ exists  in $C^{1,\beta }(\bar {\Omega})$ 
 which is then a solution for $(S)_{\lam^*}$ such that
$\mu_{1,\lambda^*}(w^*)=0.$
 In particular,   $w^*$  --often referred to as the extremal solution of problem $(S)_{\lam }$-- is unique.
\item 3.\, On the other hand, if $N\geq 8$, $f(x)=|x|^\alpha$ with $0 \le \alp \le \alp ^{**}(N):=
\frac{ 4-6N+3 \sqrt{6} (N-2) }{ 4 }$ and $\Omega$ is the unit ball,
then the extremal solution is necessarily
$w^*(x)=1-|x|^{\frac{2+\alpha}{3}}$ and is therefore singular.
\end{trivlist}}

\noindent   We remark that in general, the function $w^*$ exists in
any dimension, does solve $(S)_{\lam^*}$ in a suitable weak sense
and is the unique solution in an appropriate class.  The above
theorem says that  it is however a classical solution in dimensions
$1\leq N \leq 7$, that is
\begin{equation}\label{time:A.1}
-\Delta w^*=\frac{\lam ^*f(x)}{(1-w^*)^2}\quad  \mbox{in} \  \Omega
\,,\quad w^*>0\quad \mbox{in}\  \Omega \,,\quad w^*=0\quad
\mbox{on}\
\partial\Omega \,,
\end{equation}
and there exists an eigenfunction $\phi ^*$ of $L_{{w^*, \lam ^*}}$
satisfying
\begin{equation}\label{time:A.3}
\Delta \phi ^*+\frac{2\lam ^*\phi ^*f(x)}{(1-w^*)^3}=0\quad
\mbox{in}\  \Omega \,,\quad \phi ^*>0\quad \mbox{in}\  \Omega
\,,\quad \phi ^*=0 \quad\mbox{on}\
\partial\Omega \,.
\end{equation}
We denote by $\phi^*$ (resp., $\psi^*$) the corresponding unique
$L^2$-normalized (resp., $L^1$-normalized) positive eigenfunction of
$L_{{w^*, \lam ^*}}$.

We shall then prove in section 2 the following upper and lower
estimates on the quenching time $T=T(\lam , f,  \Omega )$ of a
solution $u$ for $(P)_\lam$ at voltage $\lambda>\lambda^*$: Under
the condition that the unique extremal solution $w^*$ of $(S)_\lam$
is regular, then

\begin{itemize}
\item For $\lam$ sufficiently close to $\lam ^*$,  we have the lower bound estimate
\begin{equation}
T(\lam ,  f, \Omega )\ge \Big(  \frac{ \sup_{x\in \Omega}\phi ^*(x) }{12\lam ^*\sup _{x\in \Omega}\frac{f(x)}{(1-w^*(x))^4}\int _\Omega \frac{\phi ^*}{(1-w^*)^2}dx}\Big)^{\frac{1}{2}} \big (\lam -\lam ^*\big)^{-\frac{1}{2}}.  \label{1:L}
\end{equation}

\item If $\int _\Omega \frac{\psi ^*(x)}{f(x)}dx<\infty$,
then
for any  $\lam >\lam ^*$,   we have the upper bound estimate
\begin{equation}
T(\lam ,  f, \Omega )\le \frac{\sqrt 3\pi }{4} \Big(\frac{\int _\Omega \frac{\psi
 ^*(x)}{f(x)}dx}{\lam
^*\int _\Omega \psi
 ^*(x)f(x)dx}\Big)^{\frac{1}{2}}\big(\lam -\lam
^*\big)^{-\frac{1}{2}}\,. \label{IV:time}
\end{equation}
\end{itemize}

\noindent Note that the above situation typically happens when
$f\equiv |x|^\beta$ and $N\le 7$, or for any $N>8$ provided $\beta$
is large. It would be interesting to establish similar estimates
 in the  case where $w^*$ is singular. In the general case, we only have the following
estimate established in section 3.
 \begin{itemize}
\item  There exist a constant $C=C(f, \Omega )>0$ and a sufficiently large $\lam _0=\lam _0 (f, \Omega )>\lam ^*$ such that for any $\lam
> \lam _0$, we have the estimates
\begin{equation}
\frac{1}{3\lam \sup_{x\in \bar \Omega }f(x)} \le T(\lam ,  f, \Omega
)\le \frac{1}{3\lam \sup_{x\in \bar \Omega }f(x)}+\frac{C}{\lam
^{\frac{2+2\alp}{2+\alp}}}\,, \label{V:time:1}
\end{equation}
where $\alp\in (0,1]$ is as in $(\ref{IV:permittivity})$.
\end{itemize}

As a byproduct of the estimate (\ref{V:time:1}),
we shall analyze and compute in section  3 that in several
situations, and at least for sufficiently large $\lam$, quenching in
finite-time must occur near the maximum point of the varying
dielectric permittivity profile $f$. More precisely, if the
quenching set  $K$ of a solution $u$ for $(P)_\lam$ is compact in
$\Omega$, and if we are in one of the following two situations:
\begin{trivlist}
\item 1)\,  $N=1$; or
\item 2)\, $N\ge 2$, $\Omega$ is a ball $B_R(0)$, $K=\{0\}$ and  $f(r)$ is radially symmetric, \end{trivlist}
then for any $a\in K$, there exists  $C>0$  such that for  $\lam $ large enough, we have
\begin{equation}
\big( \sup_{x\in\bar\Omega}f\big)^{\frac{1}{3}}-\big( f(
a)\big)^{\frac{1}{3}}\le \frac{C}{\lam ^{\frac{\alp}{2+\alp }}}\,,
\label{IV:max:1}
\end{equation}
We note that the compactness of the quenching set  has been established in \cite{GG3} (Proposition 2.1)
in the case where the domain $\Omega$ is convex and $f$ satisfies both $(\ref{IV:permittivity})$ and the additional
condition
\begin{equation}\label{permittivitya}
\begin{array}{lll}\arraycolsep=1.5pt
\hbox{$\frac{\partial f}{\partial \nu}\le 0\ $ on $\ \Omega
^c_\delta :=\{x\in \Omega :\, dist (x,\partial \Omega )\le \delta\}$
for some $\delta >0$.}
\end{array}
\end{equation}
Here $\nu$ is the outward unit norm vector to $\partial\Omega$. The
above result can be seen as a refinement of Theorem 1.1  of
\cite{GG3}   where it is proved that  under the compactness
assumption on the quenching set,  the latter set cannot contain any
zero of the  profile $f$ (see also Lemma 3.2 below).

\section{Quenching time for $\lam >\lam ^*$}

In this section, we establish the estimates on the quenching time of $(P)_\lam$. First we borrow ideas from \cite{BCM} to prove that we have quenching in  finite time as soon as $\lam >\lam ^*$, without the assumption used in \cite{GG2} that $f$ is bounded away from zero.

\begin{thm} \label{thm:2,1} 
If $\lam >\lam ^*(\Omega, f)$, then the unique solution  $u(x,t)$ of $(P)_\lam $ must quench in
finite time.
\end{thm}

\noindent{\bf Proof.} The uniqueness of solutions for $(P)_\lam$ in
$\Omega\times (0, \tau )$, where $\tau >0$ is the maximal existence
time, was already noted in Proposition 2.1 of \cite{GG2}. Let now
$\lam >\lam ^*$, and assume that  $u=u(x,t)$ of $(P)_\lam$ exists in
$\Omega \times (0,\infty )$.

Given any $0<\varepsilon <\lam -\lam ^*$, we first claim that
$(P)_{\lam -\varepsilon }$ has a global solution $u_\varepsilon$
that is uniformly bounded  in
$\Omega \times (0,\infty )$ by some constant $C_\varepsilon <1$. Indeed, set
\begin{equation}
g(u)=\frac{1}{(1-u)^2}\,,\quad h(u)=\int
^{u}_0\frac{ds}{g(s)}\,,\quad 0\le u\le 1\,, \label{finite:6}
\end{equation}
\begin{equation}
\widetilde{g}(u)=\frac{\lam -\varepsilon}{\lam (1-u)^2}\,,\quad
\widetilde{h}(u)=\int ^{u}_0\frac{ds}{\widetilde{g}(s)}\,,\quad 0\le
u\le 1\,, \label{finite:7}
\end{equation}
and let $\Phi _\varepsilon (u):=\widetilde{h}^{-1}\big(h(u)\big)$.
Direct calculations show that
\[\Phi _\varepsilon (u)=1-\big[\frac{\varepsilon}{\lam }+\frac{\lam -\varepsilon}{\lam }(1-u)^3\big]^{\frac{1}{3}}\le
C_\varepsilon<1\quad \mbox{for} \quad 0\le u\le
1\,,\] where
$C_\varepsilon =1- \big(\frac{\varepsilon}{\lam
}\big)^{\frac{1}{3}}$. Moreover, it is easy to check that $\Phi
_\varepsilon (0)=0$, that $0\le \Phi _\varepsilon (s)<s $ for $s\ge
0$, and that $\Phi _\varepsilon (s)$ is increasing and concave with
$$
\Phi _\varepsilon '(s)=\frac{\widetilde{g}(\Phi _\varepsilon(s))
}{g(s)}>0\,.
$$
Setting $v_\varepsilon =\Phi _\varepsilon (u)$, we have
$$\arraycolsep=1.5pt \begin {array}{lll} &-\Delta v_\varepsilon =-\Phi _\varepsilon ''(u)|\nabla u|^2-
\Phi _\varepsilon '(u)\Delta u&\\ [2mm] &\ge \Phi _\varepsilon
'(u)\big(\displaystyle\frac{\lam f(x)}{(1-u)^2}-u_t\big)=\lam
f(x)\Phi
_\varepsilon '(u)g(u)-(v_\varepsilon )_t &\\[2mm]
&=\lam f(x)\widetilde{g}(\Phi _\varepsilon(u))-(v_\varepsilon
)_t=\displaystyle\frac{(\lam -\varepsilon
)f(x)}{(1-v_\varepsilon)^2}-(v_\varepsilon )_t\,,&\end{array}$$ and
hence,
$v_\varepsilon =\Phi _\varepsilon (u)\le C_\varepsilon$ is therefore a
supersolution of $(P)_{\lam -\varepsilon }$. Since now zero is a
subsolution of $(P)_{\lam -\varepsilon }$, we deduce that there
exists a unique global solution $u_\varepsilon$ for $(P)_{\lam
-\varepsilon }$ satisfying  $0\le u_\varepsilon\le v_\varepsilon\le
C_\varepsilon<1$ uniformly in $\Omega \times (0,\infty )$, which
gives our first claim.

Note that $(P)_{\lam -\varepsilon }$ admits a Liapunov functional
\begin{equation} V(u_\varepsilon)=\frac{1}{2}\int _\Omega |\nabla u_\varepsilon|^2dx-(\lam
-\varepsilon )\int_\Omega \frac{f(x)}{1-u_\varepsilon}dx, \quad
\dot{V}(u_\varepsilon )=-\int _\Omega (u_\varepsilon
)^2_tdx.\label{finite:integration}
\end{equation}
Since now $\frac{1}{1-u_\varepsilon} $ is uniformly bounded in $\Omega
\times (0,\infty )$, we obtain that for $\beta <1$,
\begin{equation}
\|u_t\|_{C^{0,\beta }},\  \|u_{tt}\|_{C^{0,\beta }}<C  \quad
\mbox{uniformly bounded in}\ \Omega \times (0,\infty )\,.
\label{finite:inte}
\end{equation}
Moreover,  (\ref{finite:integration}) gives that $\int
^\infty_0\int _\Omega (u_\varepsilon )^2_tdx<\infty$, which means that $\int
_\Omega (u_\varepsilon )^2_tdx$ is a uniformly continuous function
on $[0,\infty)$, and therefore
$$
\int _\Omega (u_\varepsilon )^2_tdx\to 0\quad \mbox{as}\quad
t\to\infty\,.
$$
Further, we deduce from (\ref{finite:inte}) that $
(u_\varepsilon )_t\to 0$ as $ t\to\infty$, which shows that there
exists a function $0\le w_\varepsilon(x) <C_\varepsilon<1$ on
$\Omega$ such that $u_\varepsilon (x,t)\to w_\varepsilon (x)$ as $
t\to\infty$, where $w_\varepsilon $ satisfies
$$
-\Delta w_\varepsilon = \frac{(\lam -\varepsilon
)f(x)}{(1-w_\varepsilon )^2} \quad \mbox{in}\ \Omega, \quad
w_\varepsilon =0\quad \mbox{on}\ \partial\Omega\,.
$$
Therefore, there exists a classical solution $w_\varepsilon$ of
$(S)_{\lam -\varepsilon }$ with $\lam -\varepsilon >\lam ^*$, which
contradicts the definition of $\lam ^*$, and completes the proof of Theorem \ref{thm:2,1}.
\endproof

\subsection{Analytic estimates of  quenching time}

We now focus on estimating the quenching time $T$ when $\lam >\lam ^*$, and in the case where
the unique extremal solution $w^*$ of $(S)_\lam $ is regular. This implies that $w^*$ satisfies
\begin{equation}\label{time:2.1}
-\Delta w^*=\frac{\lam ^*f(x)}{(1-w^*)^2}\quad  \mbox{in} \  \Omega
\,,\quad w^*>0\quad \mbox{in}\  \Omega \,,\quad w^*=0\quad
\mbox{on}\
\partial\Omega \,,
\end{equation}
 and   there exists an eigenfunction
$\phi ^*$ satisfying
\begin{equation}\label{time:2.3}
\Delta \phi ^*+\frac{2\lam ^*\phi ^*f(x)}{(1-w^*)^3}=0\quad
\mbox{in}\  \Omega \,,\quad \phi ^*>0\quad \mbox{in}\  \Omega
\,,\quad \phi ^*=0 \quad\mbox{on}\
\partial\Omega \,.
\end{equation}
We shall adapt and improve some of the arguments in \cite{Lacey}.
Our first estimate is a lower bound for $T$ as stated in $(1.5)$.

\begin{thm} \label{lower:thm} Suppose that the unique extremal solution $w^*$ of $(S)_\lam$ is regular.  Then
for $\lam $ sufficiently close to $\lam ^*$,  the finite quenching time $T(\lam , f,  \Omega )$ of the unique solution $u$ for $(P)_\lam $  satisfies
\begin{equation}
T(\lam ,  f, \Omega )\ge  \Big(  \frac{ \sup_{x\in \Omega}\phi ^* }{12\lam ^*\sup _{x\in \Omega}\frac{f(x)}{(1-w^*)^4}\int _\Omega \frac{\phi ^*}{(1-w^*)^2}dx}\Big)^{\frac{1}{2}} \big (\lam -\lam ^*\big)^{-\frac{1}{2}}   \,, \label{L:L}
\end{equation}
 where $\phi ^*>0$
is the $L^2(\Omega)$-normalized eigenfunction satisfying (\ref{time:2.3}).
\end{thm}

\noindent{\bf Proof.}  Let $u^*$ be the unique solution of
$(P)_{\lam ^*}$. First, we seek a bound on the rate at which  $u^*$
approaches the corresponding steady-state $w^*$.  For that, we set
$u^*(x, t)=w^*(x)-\hat u (x, t)$. Then $\hat u(x,0)=w^*(x)$ in
$\Omega$ and $\hat u=w^*$ on $\partial\Omega$. Moreover, we have
\begin{equation}\label{L:4}
\arraycolsep=1.5pt\begin{array}{lll} \displaystyle\frac{\partial
\hat u}{\partial t}&=&\Delta \hat u-\Delta w^*-\displaystyle\frac{ \lam
 ^* f(x)}{(1-w^*+\hat u)^2}\\[3mm]
 &=&\Delta \hat u+\lam ^*f(x)\Big[\displaystyle\frac{1}{(1-w^*)^2}-\displaystyle\frac{ 1}{(1-w^*+\hat u)^2}\Big]\\[4mm]
 &\ge & \Delta \hat u+\displaystyle\frac{2\lam ^*\hat u f(x)}{(1-w^*)^3}-\displaystyle\frac{3\lam ^*\hat u^2 f(x)}{(1-w^*)^4}\\[4mm]
&\ge & \Delta \hat u+\displaystyle\frac{2\lam ^*\hat u f(x)}{(1-w^*)^3} -K_1\hat u^2 \,,
\end{array}
\end{equation}
where $K_1=3\lam ^*\sup _{x\in \Omega}\frac{ f(x)}{(1-w^*)^4}$. Define
\begin{equation}\label{L:5}
\psi =\frac{K_2\phi ^*}{t+t_0}\,,\ \quad K_2=\frac{\sup _{x\in \Omega}\phi ^*}{K_1},
 \end{equation}
where $t_0$ is chosen in such a way that
$$
\psi (x,0)=\frac{K_2\phi ^*}{t_0}\le w^*(x)=\hat u(x,0) \quad \mbox{in}\quad \Omega.
$$
Note that (\ref{L:5}) gives
$$
\Delta \psi +\displaystyle\frac{2\lam ^*\psi f(x)}{(1-w^*)^3}
-K_1\psi ^2= -\frac{K_1K_2^2}{(t+t_0)^2}(\phi ^*)^2\ge
-\frac{K_2\phi ^*}{(t+t_0)^2}=\frac{\partial \psi }{\partial t},
$$
and hence $0\le \psi \le \hat u=w^*-u^*$ in $\Omega \times (0,\infty
)$.

We now set $u=u^*+u_1$, then $u_1$ satisfies
\begin{equation}\label{L:7}
\arraycolsep=1.5pt\begin{array}{lll} \displaystyle\frac{\partial
u_1}{\partial t}&=&\Delta u_1+\displaystyle\frac{ (\lam -\lam
 ^*) f(x)}{(1- u)^2}+\lam ^*f(x)\Big[\displaystyle\frac{1}{(1-u)^2}-\displaystyle\frac{ 1}{(1-u^*)^2}\Big]\\[3mm]
 &\le &\Delta u_1+\displaystyle\frac{ (\lam -\lam
 ^*) f(x)}{(1- w^*)^2}+    \displaystyle\frac{2\lam ^*u_1f(x)}{(1-w^*)^3}   \,,
\end{array}
\end{equation}
as long as $u=u^*+u_1\le w^*$. We also define
$$
I_1=\int _\Omega \frac{\phi ^*}{(1-w^*)^2}dx\,, \quad
F(x)=\frac{f(x)}{\max \{1, \sup _{x\in \Omega}f(x)\}}\le f(x),
$$
and consider $\Phi ^*(x)\ge 0$ to be  a  nonnegative solution of the problem
\begin{equation}\label{L:8}
\arraycolsep=1.5pt\begin{array}{lll}\arraycolsep=1.5pt
   \hfill   \Delta \Phi ^* +\displaystyle\frac{2\lambda ^*
f(x)}{(1-w^*)^3}\Phi ^*+\displaystyle\frac{
f(x)}{(1-w^*)^2}-I_1\phi ^*(x)F(x)&=&0  \,   \  \ x \in \Omega, \\
 \hfill \Phi ^* (x)  &=& 0  \quad  x \in \partial\Omega\,.
\end{array}
\end{equation}
Consider also the function
\begin{equation}\label{L:10}
\psi _1=(\lam -\lam ^*)(I_1\phi ^*t+\Phi ^*) \quad \mbox{in}\quad \Omega \times (0, \tau ),
\end{equation}
where $\tau >0$ is arbitrary. Then $\psi _1(x,0)=(\lam -\lam ^*)\Phi ^*\ge 0=u_1(x,0)$ in $\Omega$, and  $\psi _1(x,t) =0=u_1(x,0)$ on $\partial \Omega$. Moreover, since $F(x)\le 1$ in $\Omega$, we obtain from $(\ref{L:7})$ and $(\ref{L:8})$ that
$$
\arraycolsep=1.5pt\begin{array}{lll} &(\psi _1-u_1)_t-\Delta  (\psi _1-u_1)&\\[2mm]
=&(\lam -\lam ^*)I_1\phi ^*-(\lam -\lam ^*)I_1t\Delta \phi ^*-(\lam -\lam ^*) \Delta \Phi ^*-(u_1)_t+\Delta u_1&\\[2mm]
\ge &(\lam -\lam ^*)I_1\phi ^*(x)-  (\lam -\lam ^*) I_1\phi ^* (x)F(x)+\displaystyle\frac{2\lambda ^*
f(x)}{(1-w^*)^3}(\psi _1-u_1) &\\[2mm]
 \ge &\displaystyle\frac{2\lambda ^*
f(x)}{(1-w^*)^3}(\psi _1-u_1) &
\end{array}
$$
in $\Omega \times (0, \tau )$, as long as $u=u^*+u_1\le w^*$. Therefore, the maximum principle implies that $\psi _1\ge u_1$ as long as $u=u^*+u_1\le w^*$.

We now obtain that
\begin{equation}\label{L:11}
u=u^*+u_1\le w^*-\psi +\psi _1= w^*-\frac{K_2\phi ^*}{t+t_0}+(\lam -\lam ^*)(I_1\phi ^*t+\Phi ^*).
\end{equation}
But the right-hand side of  $(\ref{L:11})$ is no larger than $w^*$, provided that
$$
\frac{K_2\phi ^*}{t+t_0}\ge (\lam -\lam ^*)(I_1\phi ^*t+\Phi ^*)  \quad \mbox{in}\quad \Omega \,,
$$
which is equivalent to
$$
K_2\ge (\lam -\lam ^*)(t+t_0)(I_1 t+A)\,,\quad \mbox{where}\quad  A=\sup _{x\in \Omega} \frac{\Phi ^*(x)}{\phi ^*(x)} .
$$
It requires
$$
(\lam -\lam ^*)I_1t^2+(\lam -\lam ^*)(I_1t_0+A)t-K_2+A(\lam -\lam ^*)t_0\le 0,
$$
which is
\begin{equation}\label{L:12}
t\le \displaystyle\frac{-(\lam -\lam ^*)(I_1t_0+A)+\sqrt{\Delta }}{2I_1(\lam -\lam ^*)} ,
 \end{equation}
where $$\Delta :=(\lam -\lam ^*)^2(I_1t_0+A)^2+4I_1(\lam -\lam ^*)\big(K_2-At_0(\lam -\lam ^*)\big).$$
For $\lam$  sufficiently close to $\lam ^*$,  $(\ref{L:12})$ can be satisfied if
$$t\le \frac{1}{2}\sqrt{\frac{K_2}{I_1}}(\lam -\lam ^*)^{-\frac{1}{2}}:=T_L.$$
Note that  $T_L$ is given by
$$
T_L=\Big(  \frac{ \sup_{x\in \Omega}\phi ^* (x)}{12\lam^{^*}\sup _{x\in \Omega}\frac{f(x)}{(1-w^*)^4}\int _\Omega \frac{\phi ^*}{(1-w^*)^2}dx}\Big)^{\frac{1}{2}} \big (\lam -\lam ^*\big)^{-\frac{1}{2}} .
$$
Therefore, we conclude from (\ref{L:11}) that $u\le w^*$ in $\Omega \times (0, T_L]$. This implies that the finite quenching time $T$ of $u$ satisfies $T\ge T_L$, and the proof is complete.
\endproof

We now establish the upper bound on $T$ as stated in $(1.6)$.
\begin{thm} \label{time:thm1} Suppose that the unique extremal solution $w^*$ of $(S)_\lam$ is regular,
and that $\int _\Omega \frac{\psi ^*(x)}{f(x)}dx<\infty$, where $\psi ^*>0$
is the $L^1(\Omega)$-normalized eigenfunction satisfying (\ref{time:2.3}). Then
for any  $\lam >\lam ^*$, the finite quenching time  $T=T(\lam , f,  \Omega )$ of the unique solution $u$ for $(P)_\lam $  satisfies
\begin{equation}
T(\lam ,  f, \Omega )\le \frac{\sqrt 3\pi }{4} \Big(\frac{\int _\Omega \frac{\psi
 ^*(x)}{f(x)}dx}{\lam
^*\int _\Omega \psi^*(x)
 f(x)dx}\Big)^{\frac{1}{2}}\big(\lam -\lam
^*\big)^{-\frac{1}{2}}\,. \label{IV:time:1}
\end{equation}
\end{thm}

\noindent{\bf Proof.}
Setting $u=w^*+v$, then we have
\begin{equation}\label{time:2.6}
\arraycolsep=1.5pt\begin{array}{lll} \displaystyle\frac{\partial
v}{\partial t}&=&\Delta w^*+\Delta v+\displaystyle\frac{(\lam -\lam
 ^*)f(x)}{(1-u)^2}+\displaystyle\frac{\lam
 ^*f(x)}{[1-(w^*+v)]^2}\\[4mm]
 &=&\Delta v+\displaystyle\frac{2\lam
 ^*vf(x)}{(1-w^*)^3}+\displaystyle\frac{(\lam -\lam
 ^*)f(x)}{(1-u)^2}\\[3mm]
 &&+\lam
 ^*f(x)\displaystyle\Big[\frac{1}{[1-(w^*+v)]^2}-\frac{1}{(1-w^*)^2}-\frac{2v}{(1-w^*)^3}\Big]\,.
\end{array}
\end{equation}
Multiplying $(\ref{time:2.6})$ by $\psi^*$ and integrating over
$\Omega$, we obtain
\begin{eqnarray*}\label{time:2.7}
\displaystyle\frac{d}{dt}\int _\Omega \psi^*vdx&=&(\lam -\lam
 ^*)\displaystyle\int _\Omega \frac{\psi^*f(x)}{(1-u)^2}dx\\
&&  +\lam ^*\displaystyle\int _\Omega\psi^*f(x)
 \Big[\frac{1}{[1-(w^*+v)]^2}
 -\displaystyle\frac{1}{(1-w^*)^2}-\displaystyle\frac{2v}{(1-w^*)^3}\Big]dx,
\end{eqnarray*}
where (\ref{time:2.3}) is applied.  We next define
$$
E(t)=\int _\Omega\psi^*vdx\,,\quad E(0)=-\int _\Omega\psi
^*w^*dx=-E_0\in (-1,0);
$$
$$
I_1=\int _\Omega\psi^*(x)f(x)dx\le\int _\Omega \frac{\psi
^*(x)f(x)}{(1-u)^2}dx\,,\quad I_2=\displaystyle\frac{3\lam ^*}{\int
_\Omega \frac{\psi
 ^*(x)}{f(x)}dx}.
$$
Using the inequalities
$$
\arraycolsep=1.5pt\begin{array}{lll}
\displaystyle\frac{1}{[1-(w^*+v)]^2}-\displaystyle\frac{1}{(1-w^*)^2}-\displaystyle\frac{2v}{(1-w^*)^3}
&\ge
&\left\{\arraycolsep=1.5pt\begin{array}{lll}\frac{3v^2}{(1-w^*)^4},\quad
\mbox{if} \quad v\ge 0;\\ [1mm]\frac{3v^2}{(1-u)^4},\ \, \quad
\mbox{if} \quad v\le 0;\\
\end{array}
\right.
\end{array}
$$
the H\"older inequality yields that
$$\arraycolsep=1.5pt\begin{array}{lll}
&\lam ^*\displaystyle\int _\Omega\psi^*f(x)
 \Big[\frac{1}{[1-(w^*+v)]^2}-\frac{1}{(1-w^*)^2}-\frac{2v}{(1-w^*)^3}\Big]dx&\\[4mm]
 &\ge
 3\lam ^*\displaystyle\int _\Omega v^2\psi^*(x)f(x)dx\ge \displaystyle\frac{3\lam ^*}{\int _\Omega \frac{\psi
 ^*(x)}{f(x)}dx}\Big(\int _\Omega  \psi^*vdx\Big)^2=I_2E^2(t)\,.&
\end{array}$$
It follows from the above 
that
\begin{equation}\label{time:2.8}
\frac{dE}{dt}\ge (\lam -\lam ^*)I_1+I_2E^2\,,\quad E(0)=-E_0\in
(-1,0).
\end{equation}
We now compare $E(t)$ with the solution $F(t)$ of
\begin{equation}\label{time:2.9}
\frac{dF}{dt}= (\lam -\lam ^*)I_1+I_2F^2\,,\quad F(0)=-E_0\in
(-1,0).
\end{equation}
Standard comparison principle yields that $E(t)\geq F(t)$ on their
domains of existence. Therefore,
\begin{equation}
 \sup_{\Omega} v \geq E(t) \geq  F(t) \,. \label{time:2.10}
\end{equation}
It is easy to see from (\ref{time:2.9}) that the quenching time $\bar T_1$ for $F(t)$ is
given by
\[\arraycolsep=1.5pt\begin{array}{lll}
\bar T_1 &\equiv & \Big(\displaystyle\frac{\pi }{4}+\arctan \sqrt \frac{I_2}{(\lam
-\lam ^*)I_1}\Big)\Big((\lam -\lam ^*)I_1I_2\Big)^{-\frac{1}{2}} \\ [4mm] &\le & \displaystyle
\frac{\sqrt 3\pi }{4} \Big(\displaystyle\frac{\int _\Omega \frac{\psi
 ^*(x)}{f(x)}dx}{\lam
^*\int _\Omega \psi
 ^*(x)f(x)dx}\Big)^{\frac{1}{2}}\big(\lam -\lam
^*\big)^{-\frac{1}{2}}\,.\end{array}\label{time:2.11}
\]
Therefore, for any $\lam >\lam ^*$ the unique solution $u$ of $(P)_\lam $ must
quench at a finite time $T=T(\lam , f, \Omega )\le \bar T_1 $, and we are done.
\endproof

\section{Quenching behavior  for sufficiently large $\lam$}

In this section we discuss  the quenching behavior of solutions of
$(P)_\lam $ for $\lam$ large enough. We begin with the following
refined estimates  for the quenching time as stated in $(1.7)$.

\begin{lem} \label{IV:time:lem1} Assume $f$ satisfies $(\ref{IV:permittivity})$
on a bounded domain $\Omega$, and suppose $u$ is a quenching
solution of $(P)_\lam$ at finite time $T$. Then,
there exist a constant $C=C(f, \Omega )>0$ and a sufficiently large $\lam _0=\lam _0 (f, \Omega )>0$ such that for any $\lam
> \lam _0$, we have
\begin{equation}
\frac{1}{3\lam \sup_{x\in \bar \Omega
}f(x)} \le T\le \frac{1}{3\lam \sup_{x\in \bar \Omega
}f}+\frac{C}{\lam
^{\frac{2+2\alp}{2+\alp}}}\,, \label{IV:time:1}
\end{equation}
where $\alp\in (0,1]$ is as in $(\ref{IV:permittivity})$.
\end{lem}

\noindent{\bf Proof.} In order to obtain the lower bound of finite
time $T$,  we consider the initial value problem:
\begin{equation}
\arraycolsep=1.5pt \begin {array}{lll}
\displaystyle\frac{d\eta
(t)}{dt}&=&\displaystyle\frac{\lam M}{(1-\eta (t))^2}\,, \\[2mm]
\hfill\eta (0)&=&0 \,,
  \end {array}
\label{3:14}
\end{equation}
where $M=\sup _{x\in \bar\Omega }f(x)$. From $(\ref{3:14})$ one has
$\frac{1}{\lam M}\int _0^{\eta (t)}(1-s)^2ds=t\,.$
If $T_*$ is the time where $\lim _{t\to T_*}\eta (t)=1$, then  we have
$
T_*=\frac{1}{\lam M}\int _0^{1}(1-s)^2ds=\frac{1}{3\lam M}.
$
Obviously, $\eta (t)$ is now a super-solution of $u(x,t)$ near
quenching, and thus we have
$$
T\ge T_*=\frac{1}{3\lam M}=\frac{1}{3\lam \sup _{x\in \bar\Omega
}f(x)}\,,
$$
which is true for any $\lam >0$.

We next prove the upper bound in (\ref{IV:time:1}).  Let $\bar a\in \bar\Omega $ be such that  $f(\bar
a)=\sup _{x\in\bar \Omega}f(x)$, and suppose $K=K(f,\Omega )$ is the
H\"older constant of $f$. Since $f\in C^\alpha(\bar \Omega)$
for some $\alpha \in (0,1]$, then for any sufficiently small $\varepsilon >0$, there
exists $\delta =\big(\frac{\varepsilon }{2K}\big )^{1/\alp}$ such that
$$
f(x)\ge f(\bar a)-\frac{\varepsilon}{2}\,,\quad \forall x\in
Q:=B(\bar a,\delta )\cap \Omega\,,
$$
where $B(\bar a,\delta )$ is a ball centered at $\bar a$ with radius
$\delta$. Let $v$ be the solution of
\begin{equation}
\arraycolsep=1.5pt\begin{array}{lll}\arraycolsep=1.5pt
   \hfill   v_t-\Delta v &=& \displaystyle\frac{\lambda
\big(f(\bar a)-\frac{\varepsilon}{2}\big)}{(1-v)^2} \quad \mbox{in}
\  Q\times (0, T_v)\,,
\\[2mm]
 \hfill v (x,0)  = 0\quad \mbox{in} \ &Q\,,&\quad v (x,t)  = 0\quad \mbox{on} \  \partial Q\times (0,
 T_v)\,,
\end{array}\label{IV:time:3}
\end{equation}
where $T_v$ is the maximal existence time of
$(\ref{IV:time:3})$. Comparison arguments shows that $u\ge v$ in
$Q\times (0, T_m)$, where $T_m=\min \{T,T_v\}$. Therefore, we have
$T\le T_v$.

Our goal now is to estimate $T_v$ for sufficiently large values of $\lam$. Let $\mu _1 (\delta )$ be the first eigenvalue of $-\Delta $ in
$B(\bar a,\delta )$, and let $\phi$ be the corresponding positive
eigenfunction normalized such that $\int _Q\phi dx=1$. Multiplying $(\ref{IV:time:3})$ by $\phi$ and integrating over $Q$,
we obtain
\begin{equation}
\arraycolsep=1.5pt\begin{array}{lll}\arraycolsep=1.5pt
\displaystyle\frac{d}{dt} \int_{Q} \phi v \, dx &=&
\displaystyle\int_ Q \phi \Delta v \, dx  +\displaystyle\lambda
\big(f(\bar a)-\frac{\varepsilon}{2}\big)\int_Q \frac{
\phi}{(1-v)^2} \, dx\\[3mm]
&=&-\mu _1 (\delta )\displaystyle\int_ Q \phi   v \,
dx+\displaystyle\lambda \big(f(\bar
a)-\frac{\varepsilon}{2}\big)\int_Q \frac{ \phi}{(1-v)^2} \, dx\,.
\end{array}\label{IV:time:4}
\end{equation}
Next, we define an energy-like quantity by $ E(t)=\int_Q\phi
_{_\Omega} v \, dx $ so that $E(0)=0$ and
\begin{equation}
  E(t)= \int_Q \phi_{_\Omega} v \, dx \le \sup_{Q}{v} \int_Q
 \phi  \, dx =  \sup_{Q} v \,. \label{IV:time:7}
\end{equation}
Then, using Jensen's inequality on the right-hand side of
(\ref{IV:time:4}), we obtain
\[
 \frac{dE}{dt}+\mu _1 (\delta ) E\ge \frac{\lambda \big(f(\bar
a)-\frac{\varepsilon}{2}\big) }{(1-E)^2} \,, \qquad   E(0)=0 \,.
\label{IV:time:5}
\]
Recall that there exists a constant $D=D(N)>0$, depending only on $N$,
such that $\mu _1 (\delta )=D\delta ^{-2}$. We now choose
$\varepsilon = \varepsilon (\lam , f, \Omega )
>0$ such that
\begin{equation}
\mu _1 (\delta )=D\delta ^{-2}=D\big(\frac{\varepsilon }{2K}\big )^{-\frac{2}{\alp}}=\frac{\lam}{2}
\varepsilon\,,\quad i.e.,\quad \varepsilon =\frac{2D^{\frac{\alp}{2+\alp}}K^{\frac{2}{2+\alp}}}{\lam
^{\frac{\alp}{2+\alp}}}. \label{IV:time:6}
\end{equation}
Then there exists a sufficiently large  $\lam _0=\lam _0 (f, \Omega )>\lam ^*$ such that for any $\lam
> \lam _0$, we have $f(\bar a)-\varepsilon >0$ and
$$\arraycolsep=1.5pt \displaystyle
\begin {array}{lll}
\displaystyle\frac{dE}{dt}&\ge & \displaystyle\frac{\lambda \big(f(\bar
a)-\varepsilon\big)
}{(1-E)^2}+\displaystyle\frac{\lam\varepsilon}{2(1-E)^2} -\mu _1
(\delta ) E \\[2mm ]& \ge  &\displaystyle\frac{\lambda \big(f(\bar
a)-\varepsilon\big)
}{(1-E)^2}+\displaystyle\frac{\lam\varepsilon}{2} -\mu _1 (\delta
)=\displaystyle\frac{\lambda \big(f(\bar a)-\varepsilon\big)
}{(1-E)^2}\,.\end{array}
$$
This implies a finite quenching time $T_E$ of $E$
satisfying
$$
T_E\le \frac{1}{3\lam \big(f(\bar a)-\varepsilon\big)}\le
\frac{1}{3\lam f(\bar a)}+\frac{C}{\lam
^{\frac{2+2\alp}{2+\alp}}}\,,
$$
where $C=C(f,\Omega )$ is independent of $\lam$  in view of (\ref{IV:time:6}).  Therefore, we conclude from $(\ref{IV:time:7})$ that
$$
T\le T_v\le T_E\le \frac{1}{3\lam f(\bar a)}+\frac{C}{\lam
^{\frac{2+2\alp}{2+\alp}}}\,,
$$and the lemma is proved.
\endproof

We now recall the following result proved in Theorem 1.1 of
\cite{GG3}.

\begin{lem} \label{lem:rate1}
 Assume $f$ satisfies
$(\ref{IV:permittivity})$ for some $\alp\in (0,1]$ on a bounded
domain $\Omega\subset\mathbb{R}^N$, and let $u$ be a quenching solution of
$(P)_\lam$ at finite time $T$. Assuming the quenching set of $u$ is compact in $\Omega$, then
\begin{enumerate}
\item  No point $a\in \bar\Omega$ satisfying $f(a)=0$ can be a quenching point of $u$;
\item There exists a  constant
$M>0$ such that
\begin{equation}
M(T-t)^{\frac{1}{3}}\le 1-u(x,t)\quad \mbox{in}\quad \Omega\times
(0,T).\label{rate:upper1.1}
\end{equation}
\end{enumerate}
\end{lem}

The following result can now be seen as a converse of Lemma 3.2:
for sufficiently large $\lam$, finite-time quenching  must occur
near the maximum point of the varying dielectric permittivity
profile $f$.

\begin{thm} \label{IV:time:thm1}  Assume $f$ satisfies
$(\ref{IV:permittivity})$ for some $\alp\in (0,1]$ on a bounded
domain $\Omega\subset \mathbb{R}^N$, and suppose that $u$ is a quenching solution of
$(P)_\lam$ at finite time $T$, in such a way that the quenching set  $K$ of $u$ is compact in  $\Omega$. Then, for any $a\in K$, there exists  $C>0$  such that for  $\lam $ large enough, we have
\begin{equation}
\big( \sup_{x\in\bar\Omega}f\big)^{\frac{1}{3}}-\big( f(
a)\big)^{\frac{1}{3}}\le \frac{C}{\lam ^{\frac{\alp}{2+\alp }}}\,,
\label{IV:max:1}
\end{equation}
provided we are in one the following two situations:
\begin{trivlist}
\item 1)\,  $N=1$; or
\item 2)\, $N\ge 2$ and $a=0$, $\Omega$ is a ball $B_R(0)$ and $f(r)$ is radially symmetric. \end{trivlist}
\end{thm}

\noindent{\bf Proof.} The idea of the proof --inspired by
\cite{CER}-- is to combine the estimates on quenching time given by
Lemma \ref{IV:time:lem1}, with the local energy estimates near any
quenching point established in \cite{GG3}.  Given a quenching point
$a$ of $u$ and its corresponding quenching time $T$,
we define
$$
y=\frac{x-a}{\sqrt{T-t}}\,, \quad s=-\log(1-\frac{t}{T})\,,\quad
1-u(x,t)=(T-t)^{\frac{1}{3}}w(y,s),
$$
then $w$ satisfies
$$
\rho w_s=\nabla \cdot (\rho \nabla w)+\frac{1}{3}\rho w-\frac{\lam
\rho f(a+yT^{\frac{1}{2}}e^{-\frac{s}{2}})}{w^2}\quad \mbox{in}\quad
\Omega (s)\times (0,\infty )\,,
 $$
where $\rho (y)=e^{-|y|^2/4}$ and $\Omega (s)=\{y:
a+yT^{\frac{1}{2}}e^{-\frac{s}{2}}\in \Omega\}$. The compactness assumption on the quenching set implies that there
exists a sufficiently large $s_0>0$ such that $B_s(a)\subset \Omega (s)$ for any $s\ge
s_0$.

Consider now the ``frozen" energy functional
$$
E(w)= \frac{1}{2}\int _{B_s} \rho |\nabla w|^2dy-\frac{1}{6}\int
_{B_s} \rho  w^2dy- \int_{B_s} \frac{\lam \rho f(a)}{w}dy\,,
$$ which is defined in the compact set $B_s$ of $\Omega
_a(s)$ for $s\ge s_0$. Note from Lemma \ref{lem:rate1} that $f(a)>0$.  Using the same argument of Lemma 2.10 in
\cite{GG3}, one can obtain
\begin{equation}\arraycolsep=1.5pt\begin{array}{lll}\arraycolsep=1.5pt
\displaystyle\int _{B_s}\rho |w_s|^2dy &\le
&-\displaystyle\frac{dE}{d s}+\displaystyle\int _{\partial
B_s}\rho w_s\frac{\partial w}{\partial
\nu}dS+\displaystyle\frac{1}{2s}\int _{\partial B_s}\rho
|\nabla w|^2(y\cdot \nu )dS\\[4mm]
&&+\displaystyle\int _{B_s}\frac{\lam\rho
w_s[f(a)-f(a+yT^{\frac{1}{2}}e^{-\frac{s}{2}})]}{w^2}dy\\[4mm]
&:=&-\displaystyle\frac{dE}{d s}+I_1+I_2+I_3\,,
\end{array}\label{IV:time:8}
\end{equation}
where $$ I_1\le C_1s^Ne^{-\frac{s^2}{4}+\frac{s}{3}}\,,\quad I_2 \le
C_3s^{N-1}e^{-\frac{s^2}{4}}.
$$
To estimate $I_3$, we use Lemma \ref{lem:rate1} to infer that $w$ has a lower bound, and since $f\in
C^\alp (\bar\Omega )$, we apply H\"older's
inequality to deduce that
$$
I_3\le CT^{\frac{\alp}{2}}e^{-\frac{\alp}{2}s}\int _{B_s}\rho
|y|^\alp w_sdy\le CT^{\frac{\alp}{2}}e^{-\frac{\alp}{2}s}\Big(\int
_{B_s}\rho |w_s|^2dy\Big)^{\frac{1}{2}}.
$$
Therefore, (\ref{IV:time:8}) gives for $s\gg 1$,
\begin{equation}
\displaystyle\frac{dE}{d s}\le -\displaystyle\int _{B_s}\rho |w_s|^2dy+CT^{\frac{\alp}{2}}e^{-\frac{\alp}{2}s}\Big(\int
_{B_s}\rho |w_s|^2dy\Big)^{\frac{1}{2}}+Cs^{N}e^{-\frac{s^2}{4}+\frac{s}{3}}.
\label{IV:time:8'}
\end{equation}
Maximizing now the right hand side of (\ref{IV:time:8'}) with respect to $\displaystyle\int _{B_s}\rho |w_s|^2dy$, it yields that for  $s\gg 1$
$$\frac{dE}{ds}\le CT^\alp e^{-\alp s} +Cs^{N}e^{-\frac{s^2}{4}+\frac{s}{3}}\le CT^\alp e^{-\alp s}.$$
This leads to
$$E(w)\le E\big(w(y,0)\big)+\frac{CT^\alp}{\alp} =E(T^{-\frac{1}{3}})+\frac{CT^\alp}{\alp} \,.$$
Under the compactness assumption on the quenching set, a proof similar to Theorem 1.3 in \cite{GG3} (see also \cite{GK,GK2})  gives
that
$$
\lim _{s\to\infty}w(y,s)= \big(3\lam
f(a)\big)^{\frac{1}{3}} :=k(a)
$$
uniformly on $|y|\le C$ for any bounded constant $C$, and $E(w(\cdot
,s))\to E(k(a))$ as $s\to \infty$, provided  one of the following
conditions holds:
\begin{trivlist}
\item 1)\,  $N=1$; or
\item 2)\, $N\ge 2$ and $a=0$, $\Omega =B_R(0)$ is a bounded
ball and  $f(r)=f(|x|)$ is radially symmetric. \end{trivlist}
Therefore, under the assumption of Theorem \ref{IV:time:thm1},  we have the following upper bound
\begin{equation}
E(k(a))\le E(T^{-\frac{1}{3}})+\frac{CT^\alp}{\alp} \,. \label{IV:time:9}
\end{equation}

Observe that if $b$ is a constant then the energy $E$ can be
rewritten as $E(b)=\Gamma F(b)$, where $\Gamma =\int \rho (y) dy$
and $F$ is the function
$$F(z)=-\frac{1}{6}z^2-\frac{\lam f(a)}{z}\,,\quad z>0\,.$$
Since $F$ attains a unique maximum at $k(a)$ and $F''(k(a))=-1$, there
exist $\gamma$ and $\beta $  such that if $|z-k(a)|\le \gamma $ then
$F''(z)\le -\frac{1}{2}$, and if $|F(z)-F(k(a))|\le \beta$ then
$|z-k(a)|\le \gamma$. So we obtain from $(\ref{IV:time:9})$ that
$$
F(k(a))\le F(T^{-\frac{1}{3}})+\frac{CT^\alp}{\alp}\,.
$$
Choose $\lam _1$ such that $\frac{CT^\alp}{\alp} =\beta$. Then for
$\lam > \max\{\lam _0,\lam _1\}$, where $\lam _0$ is as in Lemma
\ref{IV:time:lem1}, we have
$$
\beta \ge\frac{CT^\alp}{\alp}\ge F(k(a))-F(T^{-\frac{1}{3}})\,.
$$
Hence from the properties of $F$, we have $k(a)-T^{-\frac{1}{3}}\le
\gamma$, which implies $F''(k(a))\le -\frac{1}{2}$. It now deduces from (\ref{IV:time:9}) that
$$
\frac{C}{\alp \lam ^\alp}\ge \frac{CT^\alp}{\alp}\ge F(k(a))-F(T^{-\frac{1}{3}})\ge \frac{1}{4}[T^{-\frac{1}{3}}-k(a)]^2,
$$ where Lemma \ref{IV:time:lem1} is applied in the first inequality. This further gives that
\begin{equation}
T^{-\frac{1}{3}}-\big(3\lam
f(a)\big)^{\frac{1}{3}} \le \frac{C}{\lam ^{\frac{\alp}{2}}}\,. \label{IV:time:10}
\end{equation}
On the other hand, since  Lemma \ref{IV:time:lem1} gives
$$
T\le \frac{1}{3\lam \sup _{x\in\bar \Omega}f}+\frac{C}{\lam
^{\frac{2+2\alp}{2+\alp}}}\le \frac{1}{3\lam \sup _{x\in\bar
\Omega}f}\Big(1+\frac{C}{\lam ^{\frac{\alp}{2+\alp}}}\Big)\,,
$$
we have
$$
T^{-\frac{1}{3}}\ge \big(3\lam \sup _{x\in\bar \Omega}f(x)\big)^{\frac{1}{3}}\Big(1-\frac{C}{\lam
^{\frac{\alp}{2+\alp}}}\Big)\,.
$$
Therefore, we finally conclude  that
$$
\big(  \sup _{x\in\bar \Omega}f(x)\big)^{\frac{1}{3}}-\big(  f(
a)\big)^{\frac{1}{3}}\le \frac{C}{\lam ^{\frac{1}{3}+\frac{\alp}{2}}}+\frac{C}{\lam
^{\frac{\alp}{2+\alp}}}\le \frac{C}{\lam ^{\frac{\alp}{2+\alp}}}.
$$
This completes the proof of Theorem \ref{IV:time:thm1}.
\endproof

\begin{figure}[htbp]
\begin{center}
{\includegraphics[width = 7.88cm,height=4.6cm,clip]{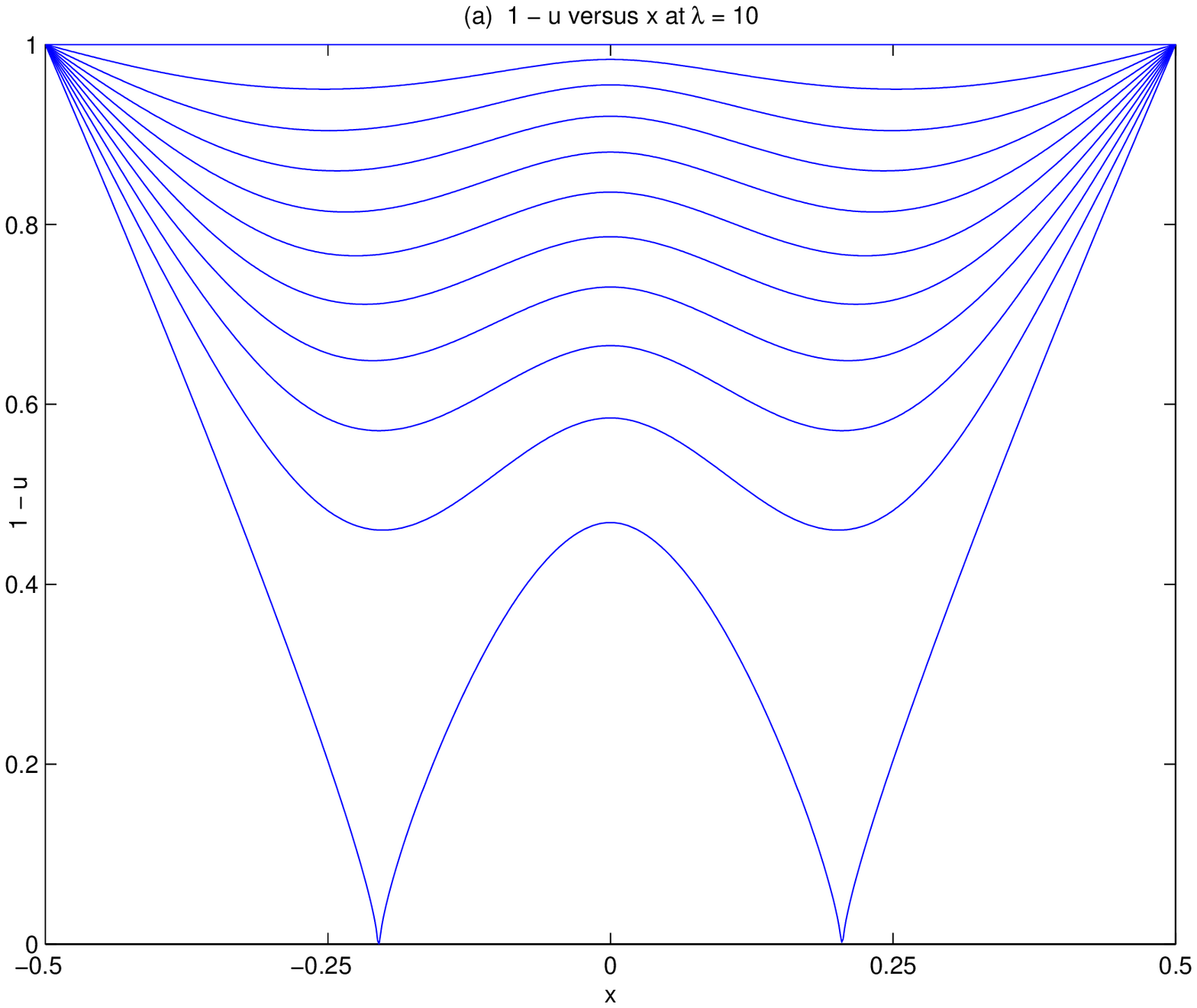}
\label{fig:fig2_2a}}
{\includegraphics[width = 7.88cm,height=4.6cm,clip]{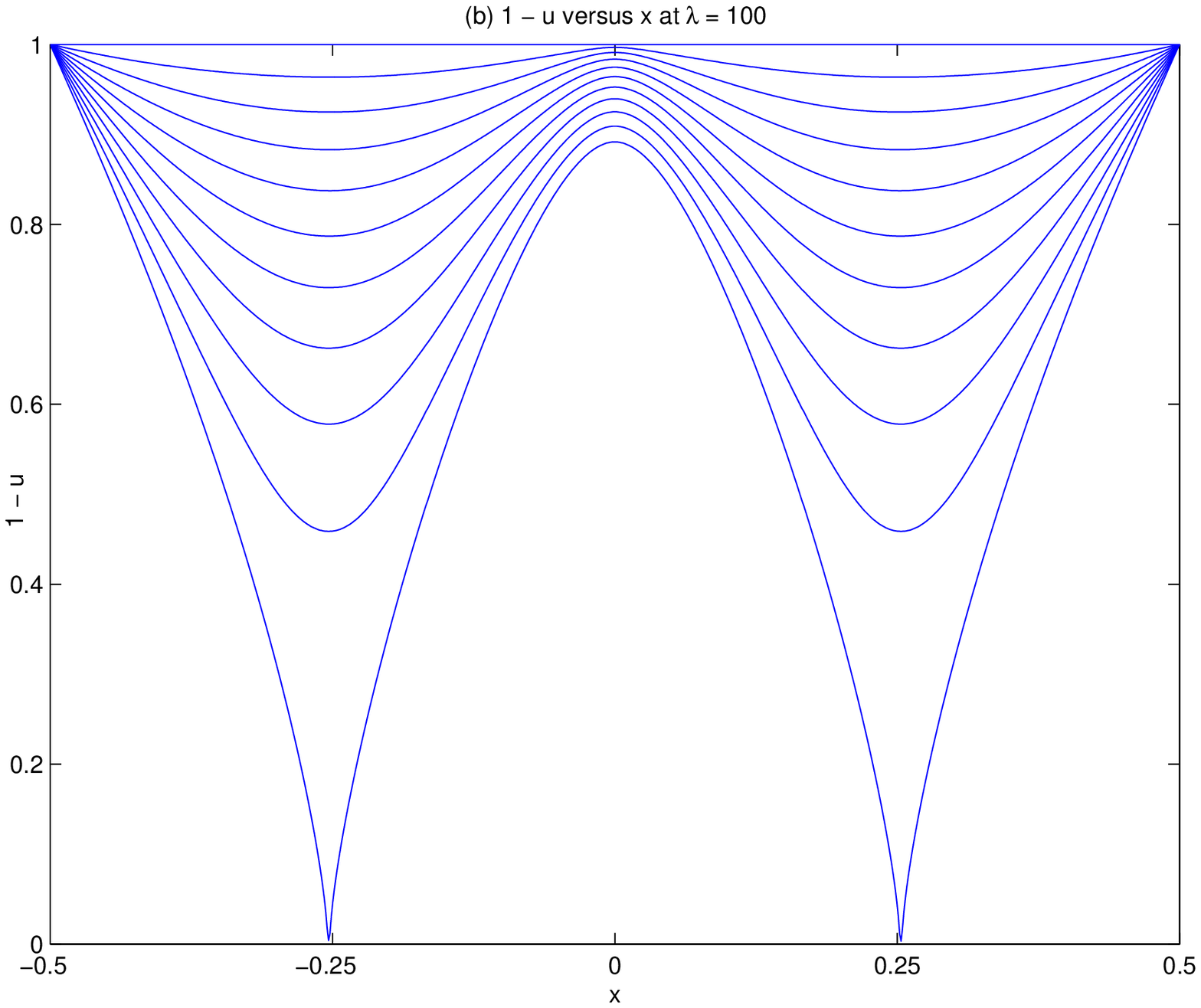}
\label{fig:fig2_2b}} \caption{{\em Left figure $(a)$:   plots of
$1-u$ versus $x$ at different times, where  $\lam =10$. Right figure
$(b)$: plots of $1-u$ versus $x$ at different times, where  $\lam
=100$.  }} \label{fig:fig2_2}
\end{center}
\end{figure}

Before ending this section, we now present a few numerical
simulations on Lemma 3.1 and Theorem \ref{IV:time:thm1}. Here we
apply the implicit Crank-Nicholson scheme (see \S 3.2 of \cite{GPW}
for details), with the meshpoints $N=6000$, to $(P)_\lam$ in the
symmetric slab domain $-1/2\le x\le 1/2$. We choose the varying
dielectric permittivity profile $f(x)$ satisfying
\begin{equation}
f[\alpha ](x)=\left \{ \begin {array}{ll}
 1-16(x+1/4)^2\,, \ \ \ & \mbox{if} \ \ \  x< -1/4 \,;  \\
 |\sin(2\pi x)|\,, \ \ \ & \mbox{if} \ \ \  |x|\le 1/4 \,;  \\
  1-16(x-1/4)^2\,, \ \ \ & \mbox{if} \ \ \ \, x> 1/4 \,.
 \end {array}
 \right.
\label{5:3}
\end{equation}
Note that $x=\pm 0.25$ are two maximum points of $f(x)$, and all
assumptions of Lemma 3.1 and Theorem \ref{IV:time:thm1} are
satisfied in view of $(1.9)$.

\begin{figure}[htbp]
\begin{center}
{\includegraphics[width = 7.88cm,height=4.6cm,clip]{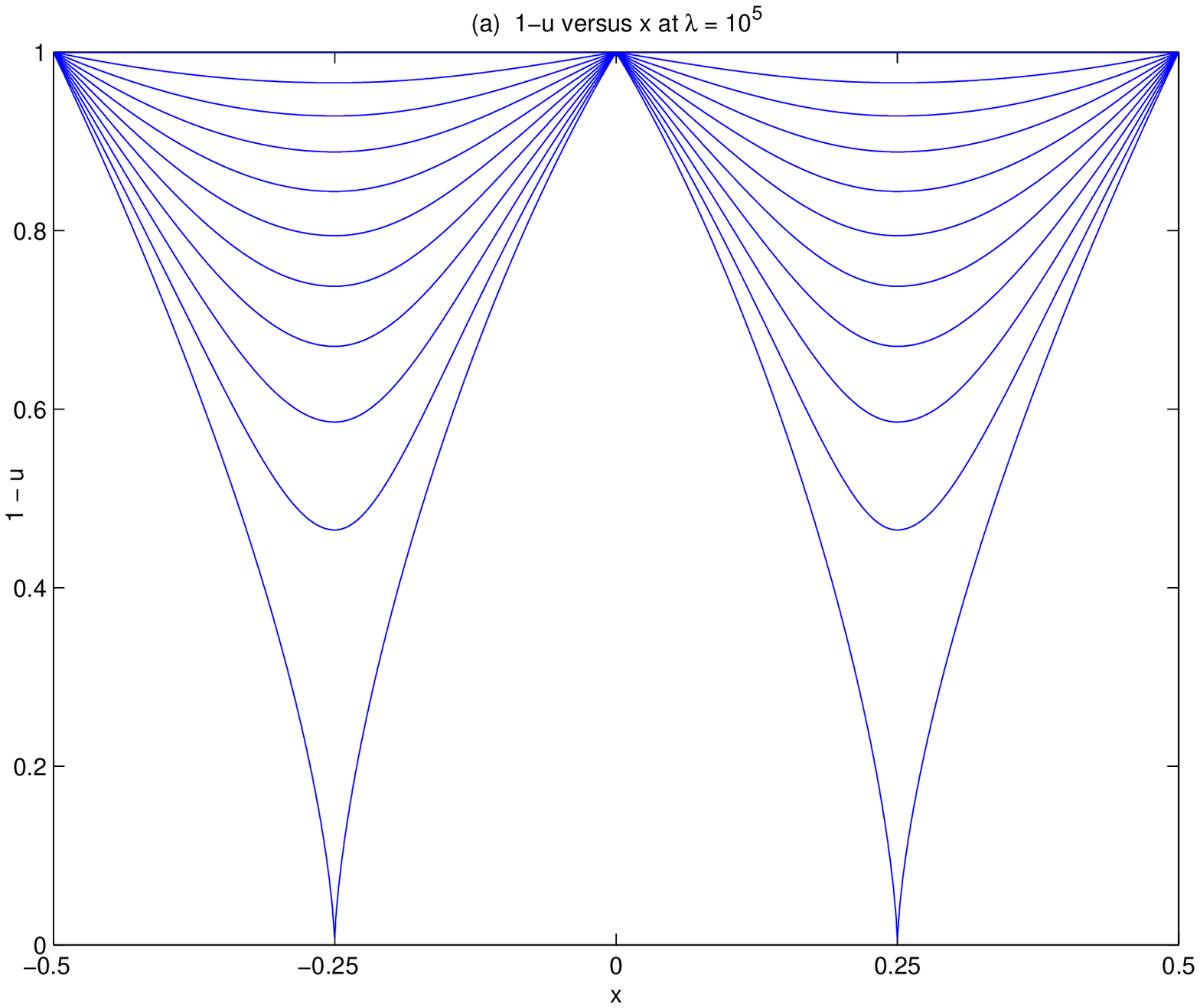}
\label{fig:fig2_2a}}
{\includegraphics[width = 7.88cm,height=4.6cm,clip]{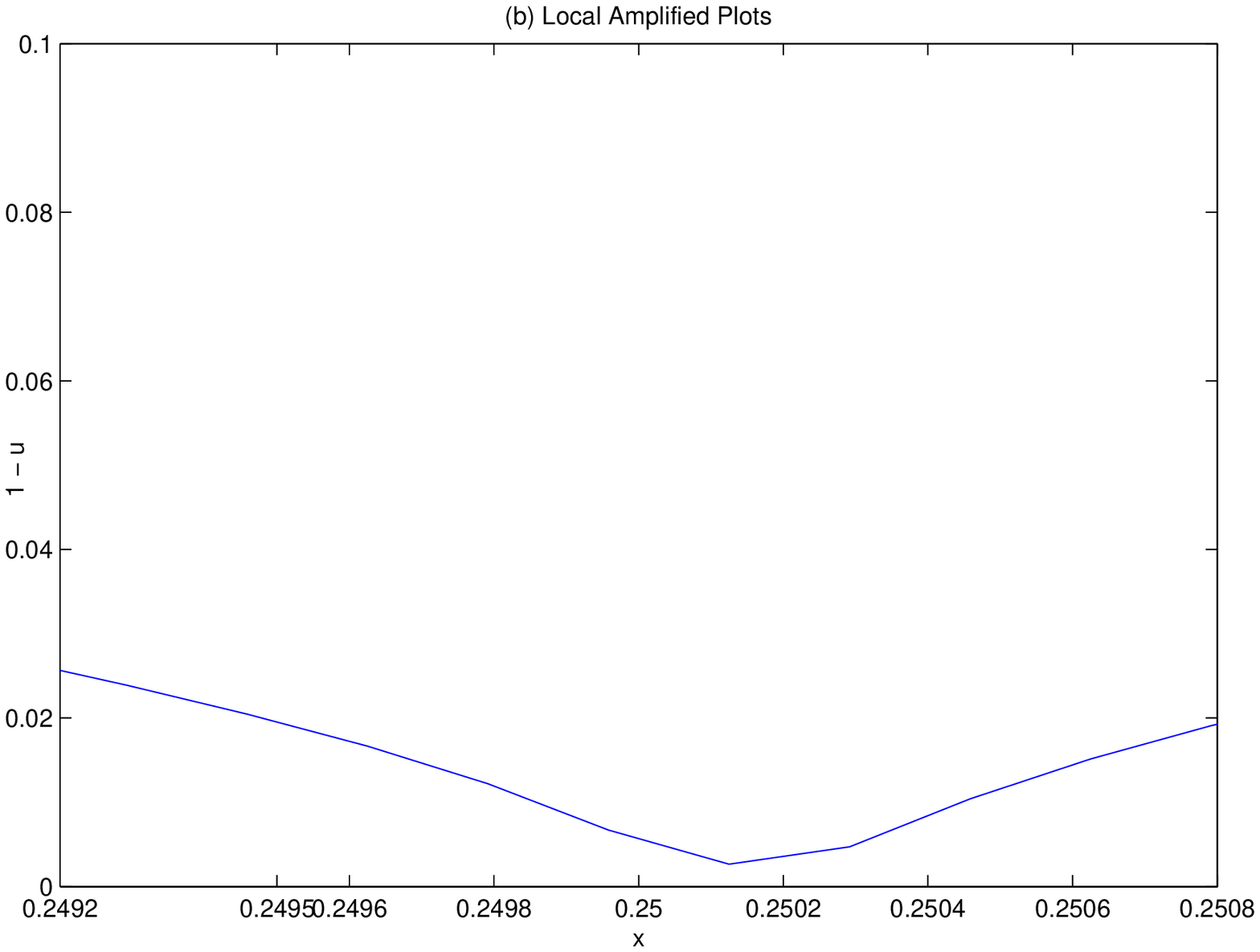}
\label{fig:fig2_2b}} \caption{{\em Left figure $(a)$: plots of $1-u$
versus $x$ at different times. Right figure $(b)$: local amplified
plots of $(a)$.   }} \label{fig:fig2_2}
\end{center}
\end{figure}


\noindent {\bf Simulation 1.} Quenching behavior for small   $\lam> \lam ^*$:\\
In Fig. 2$(a)$: $1-u$ versus $x$ is plotted at different times for
$(P)_\lam$ at $\lambda =10$, where the quenching time is $ T=
0.05174132$. The quenching   is observed at $x= \pm 0.204 $, a bit
far away from the maximum points of profile $f(x)$. In Fig.~2$(b)$:
$1-u$ versus $x$ is plotted at different times for $(P)_\lam$ at
$\lambda =100$, where the quenching time is $ T= 0.003523908$. In
this case, the quenching   is observed at  $x= \pm  0.2535 $, very
close to the maximum points of profile $f(x)$. This simulation shows
the necessary of the assumption that   Lemma 3.1 and  Theorem
\ref{IV:time:thm1}  hold only for  sufficiently large $\lam$.

\noindent {\bf Simulation 2:} Quenching behavior for sufficiently large $\lam$:\\
In Fig. 3$(a)$, $1-u$ versus $x$ is plotted at different times for
$(P)_\lam$ at $\lambda =10^5$, where   the quenching time is $ T=
0.000003332783$.  In this case, two quenching points are observed at
$x= \pm 0.250165 $, more close to the maximum points of profile
$f(x)$. In Fig.~3$(b)$ we show the local amplified plots of $(a)$
near the maximum point $x=0.25$ of $f(x)$. By further increasing the
value of $\lam $, we observe that quenching points become further
close to the maximum points of $f(x)$.

\end{document}